\documentclass[a4paper,11pt]{article}
\pdfoutput=1
\usepackage{amsmath, amsthm, amssymb}
\usepackage{graphicx}
\usepackage{color}
\usepackage{epsfig}
\usepackage{a4wide}
\usepackage{verbatim}
\newtheorem{theorem}{Theorem}

\newtheorem{lemma}[theorem]{Lemma}

\newtheorem{proposition}[theorem]{Proposition}
\newtheorem{remark}[theorem]{Remark}

\begin{document}
\title{Further calculations for the McKean stochastic game for a spectrally negative L\'evy process: from a point to an interval}
\author{\textbf{E.J. Baurdoux\footnote{Department of Statistics, London School of Economics. Houghton street, {\sc London, WC2A 2AE, United Kingdom.} E-mail: e.j.baurdoux@lse.ac.uk}, K. van Schaik\footnote{Department of Mathematical Sciences, University of Bath, Claverton Down, {\sc Bath, BA2 7AY, United Kingdom}. E-mail: k.van.schaik@bath.ac.uk. This author gratefully acknowledges being supported by a post-doctoral grant from the AXA Research Fund}}}

\date{}

\maketitle

\begin{abstract}
Following Baurdoux and Kyprianou \cite{McKean} we consider the McKean stochastic game, a game version of the McKean optimal stopping problem (American put), driven by a spectrally negative L\'evy process. We improve their characterisation of a saddle point for this game when the driving process has a Gaussian component and negative jumps. In particular we show that the exercise region of the minimiser consists of a singleton when the penalty parameter is larger than some threshold and `thickens' to a full interval when the penalty parameter drops below this threshold. Expressions in terms of scale functions for the general case and in terms of polynomials for a specific jump-diffusion case are provided.
\end{abstract}

\vspace{0.1cm}
\noindent
{\footnotesize Keywords:} {\footnotesize Stochastic games, optimal stopping, L{\'e}vy processes, fluctuation theory}

\vspace{0.1cm}
\noindent
{\footnotesize Mathematics Subject Classification (2000): 60G40, 91A15}

\section{Introduction}

This paper is a follow-up to the paper \cite{McKean} by Baurdoux and Kyprianou (henceforth BK), in which the solution to the McKean stochastic game driven by a spectrally negative L\'evy process is studied.
Let us introduce the setting in BK (and in this paper). Let $X$ be a L\'evy process defined on a filtered
probability space $(\Omega,\mathcal{F},\mathbf{F},\mathbb{P})$, where $\mathbf{F}=(\mathcal{F}_t)_{t \geq 0}$ is the filtration generated by $X$ which is naturally enlarged (cf. Definition 1.3.38 in Bichteler \cite{Bichteler02}). For $x \in \mathbb{R}$ we denote by $\mathbb{P}_x$ the law of $X$ when it is
started at $x$ and we abbreviate $\mathbb{P}=\mathbb{P}_0$. Accordingly we shall write $\mathbb{E}_x$ and $\mathbb{E}$ for the
associated expectation operators. We assume throughout that $X$
is spectrally negative, meaning that it has no positive jumps and that it is not
the negative of a subordinator.

The McKean stochastic game is an example of a type of stochastic games introduced by Dynkin \cite{Dynkin69}. It is a two-player zero sum game, consisting of a maximiser aiming at maximising over $\mathbf{F}$-stopping times $\tau$ the expected payoff according to the (discounted) lower payoff process given by $e^{-qt}(K-\exp(X_t))^+$ for all $t \geq 0$ and a minimiser aiming at minimising over $\mathbf{F}$-stopping times $\sigma$ the expected payoff according to the (discounted) upper payoff process given by $e^{-qt}((K-\exp(X_t))^+ +\delta)$ for all $t \geq 0$, where $K,\delta>0$. That is, for any pair of stopping times $(\tau,\sigma)$ the payoff to the maximiser is

\[ M_x(\tau,\sigma) := \mathbb{E}_x [ e^{-q \tau} (K-e^{X_{\tau}})^+ \mathbf{1}_{\{ \tau \leq \sigma \}} + e^{-q \sigma} ((K-e^{X_{\sigma}})^+ +\delta) \mathbf{1}_{\{ \sigma < \tau \}} ]. \]
We assume throughout this paper that the discount factor $q$ satisfies

\begin{equation}
0\leq \psi(1)\leq q \text{ and } q>0,
\label{ass}
\end{equation}
where $\psi$ denotes the Laplace exponent of $X$. For a spectrally negative L\'evy process $X$ this Laplace exponent is of the form
\[\psi(\lambda)=a\lambda+\frac{1}{2}\sigma_X^2\lambda^2+\int_{(-\infty,0)}(e^{\lambda x}-1-\lambda x 1_{\{x>-1\}})\Pi(dx)\]
for $\lambda\geq 0$.
Here $\sigma_X$ is called the Gaussian coefficient of $X$ and $\Pi$ is the L\'evy measure which satisfies $\int_{(-\infty,0)}(1\wedge x^2)\Pi(dx)<\infty$.

 Note that since both payoff processes vanish a.s. as $t \to \infty$, there is no ambiguity in allowing for $\tau$ and $\sigma$ to be infinitely valued as we will in this paper. For any $x$, this game has a \emph{value} if the upper and lower value, $\inf_{\sigma} \sup_{\tau} M_x(\tau,\sigma)$ and $\sup_{\tau}  \inf_{\sigma} M_x(\tau,\sigma)$ respectively, coincide. Even more, if a pair $(\tau^*,\sigma^*)$ exists such that

\[ M_x(\tau,\sigma^*) \leq M_x(\tau^*,\sigma^*) \leq M_x(\tau^*,\sigma) \quad \mbox{for all $(\tau,\sigma)$}, \]
the value exists and equals $M_x(\tau^*,\sigma^*)$. In this case $(\tau^*,\sigma^*)$ is called a \emph{saddle point} (or Nash equilibrium). For an account of these concepts in a general Markovian setting, see Ekstr\"om and Peskir \cite{Ekstrom08} and the references therein. For other examples of stochastic games, see e.g. Kifer \cite{Kifer00}, Kyprianou \cite{Kyprianou04}, Baurdoux and Kyprianou \cite{Baurdoux08}, Gapeev and K\"uhn \cite{Gapeev05}, Baurdoux et al. \cite{Baurdoux09}.

Note that the McKean game can be seen as an extension of the classic McKean optimal stopping problem (cf. \cite{McKean65} and Theorem \ref{avram} below). In a financial interpretation, this optimal stopping problem is usually referred to as American put option, with $K$ the strike price. The McKean game then extends the American put option by introducing the possibility for the writer of the option to cancel the contract, at the expense of paying the intrinsic value plus an extra constant penalty given by the penalty parameter $\delta$. Cf. e.g. Kifer \cite{Kifer00} and Kallsen and K\"uhn \cite{Kallsen04} for a general account on the interpretation of stochastic games as financial contracts.

In BK it was shown that a saddle point $(\tau^*,\sigma^*)$ indeed exists for the McKean game, so in particular the value function $V$ is well defined by

\begin{eqnarray}
V(x) &=& \sup_\tau \inf_\sigma \mathbb{E}_x \left[ e^{-q\tau}(K-e^{X_\tau})^+\mathbf{1}_{\{\tau \leq \sigma\}}+e^{-q\sigma}((K-e^{X_\sigma})^++\delta)\mathbf{1}_{\{\sigma<\tau\}} \right] \nonumber\\
&=& \inf_\sigma \sup_\tau \mathbb{E}_x \left[ e^{-q\tau}(K-e^{X_\tau})^+\mathbf{1}_{\{\tau \leq \sigma\}}+e^{-q\sigma}((K-e^{X_\sigma})^++\delta)\mathbf{1}_{\{\sigma<\tau\}} \right] \nonumber\\
&=& \mathbb{E}_x \left[ e^{-q\tau^*}(K-e^{X_{\tau^*}})^+\mathbf{1}_{\{\tau^* \leq \sigma^*\}}+e^{-q\sigma^*}((K-e^{X_{\sigma^*}})^++\delta)\mathbf{1}_{\{\sigma^*<\tau^*\}} \right]. \nonumber
\end{eqnarray}
The optimal stopping time for the maximiser, $\tau^*$, is the first hitting time of an interval of the form $(-\infty,x^*]$ for some $x^*<\log K$. For the minimiser the optimal stopping time $\sigma^*$ is as follows. When the penalty parameter $\delta$ exceeds $\bar{\delta}:=U(\log K)$, where $U$ denotes the value function of the McKean optimal stopping problem, the minimiser never stops (i.e. $\sigma^*=\infty$). When $\delta < \bar{\delta}$, the optimal stopping region for the minimiser is an interval of the form $[\log K,y^*]$. If the Gaussian component $\sigma_X$ of $X$ is equal to zero (note that this corresponds to the situation that $X$ does not creep downwards), we have $y^*>\log K$. Furthermore formulae in terms of scale functions for $x^*$ and $V$ on $(-\infty,\log K]$ were provided.

However, two issues were left open in BK. Firstly, when $X$ has a Gaussian component it was not clear when the optimal stopping region for the minimiser consists of a point and when of an interval, i.e. when $y^*=\log K$ and when $y^*>\log K$ holds. Secondly, no characterisation was given of $y^*$. In this paper we give an answer to both these issues. In particular, we show that when $\sigma_{X}>0$ there exists a critical value $\delta_0 \in (0,\bar{\delta})$ such that the stopping region for the minimiser is a single point when $\delta\in[\delta_0,\bar{\delta})$ and a full interval when $\delta\in(0,\delta_0)$, cf. Theorem \ref{thm_Erik_main} (see also Remark \ref{rem_Kees1}). Furthermore we show that $y^*$ and $\delta_0$ can be characterised as unique solutions to functional equations using scale functions, cf. Theorem \ref{thm_Erik2}.

The rest of this paper is organised as follows. In the remainder of this introduction we introduce scale functions and some notation (Subsection \ref{subsec_scale}), and review the results from BK in more detail (Subsection \ref{subsec_review}). In Section \ref{sec_main_res} we present our new results. Finally, in Section \ref{sec_jump_diff} we translate these results to a specific jump-diffusion setting, accompanied by some plots.

\subsection{Scale functions}\label{subsec_scale}

First we introduce some notation for first entry times. For $a \leq b$ we write

\[ \tau_a^+ := \inf \{ t>0 \, | \, X_t>a \}, \quad \tau_a^- := \inf \{ t>0 \, | \, X_t<a \} \quad \mbox{and} \quad T_{[a,b]} := \inf \{ t>0 \, | \, X_t \in [a,b] \}. \]
Furthermore we denote the often used first hitting time of $\log K$ for simplicity by $T_K$, that is $T_K := \inf \{ t>0 \, | \, X_t = \log K \}$.

A useful class of functions when studying first exit problems driven by spectrally negative L\'evy processes are so-called scale functions. We shortly review some of their properties as they play an important role in this paper, for a more complete overview the reader is e.g. referred to Chapter VII in Bertoin \cite{Bertoin96} or Chapter 8 in Kyprianou \cite{Kyprianou06}.
For each $q \geq 0$ the scale functions $W^{(q)}: \mathbb{R} \to [0,\infty)$ are known to satisfy for all $x \in \mathbb{R}$ and $a \geq 0$

\begin{equation}\label{K_30sept2}
 \mathbb{E}_x [ e^{-q \tau_a^+} \mathbf{1}_{\{ \tau_a^+<\tau_0^- \}} ] = \frac{W^{(q)}(x \wedge a)}{W^{(q)}(a)}.
\end{equation}
In particular it is evident that $W^{(q)}(x)=0$ for all $x<0$. Furthermore it is known that $W^{(q)}$ is almost everywhere differentiable on $(0,\infty)$, it is right continuous at zero and

\begin{equation}\label{K_18sept2}
\int_0^{\infty} e^{-\beta x} W^{(q)}(x) \, \mbox{d}x = \frac{1}{\psi(\beta)-q}
\end{equation}
for all $\beta>\Phi(q)$, where $\Phi(q)$ is the largest root of the equation $\psi(\theta)=q$ (of which there are at most two, recall that $\psi$ is the Laplace exponent of $X$). We shall assume throughout this paper that the jump measure $\Pi$ has no atoms when $X$ is of bounded variation, which implies that $W^{(q)}\in C^1(0,\infty)$ (see \cite{Doney}).
In case $X$ has a Gaussian component $\sigma_X>0$ it is known that $W^{(q)} \in C^2(0,\infty)$ with $W^{(q)}(0)=0$ and $W^{(q){\prime}}(0)=2/\sigma_X^2$. We usually write $W=W^{(0)}$.

Associated to the functions $W^{(q)}$ are the functions $Z^{(q)}: \mathbb{R} \to [1,\infty)$ defined by

\begin{equation}\label{K_18sept3}
Z^{(q)}(x) = 1+q \int_0^x W^{(q)}(y) \, \mbox{d}y
\end{equation}
for $q \geq 0$. Together the functions $W^{(q)}$ and $Z^{(q)}$ are collectively known as scale functions and predominantly appear in almost all fluctuation identities for spectrally negative L\'evy processes. For example, it is also known that for all $x \in \mathbb{R}$ and $a,q \geq 0$

\[ \mathbb{E}_x [ e^{-q \tau_0^-} \mathbf{1}_{\{ \tau_a^+>\tau_0^- \}} ] = Z^{(q)}(x \wedge a) - \frac{Z^{(q)}(a)}{W^{(q)}(a)} W^{(q)}(x \wedge a) \]
and

\begin{equation}\label{K_30sept1}
\mathbb{E}_x [ e^{-q \tau_0^-} \mathbf{1}_{\{ \tau_0^-<\infty \}} ] = Z^{(q)}(x) - \frac{q}{\Phi(q)} W^{(q)}(x),
\end{equation}
where $q/\Phi(q)$ is to be understood in the limiting sense $\psi'(0) \wedge 0$ when $q=0$.

For $c > 0$, consider the change of measure

\begin{equation}\label{K_18sept1}
\left. \frac{\mathrm{d}\mathbb{P}^c}{\mathrm{d}\mathbb{P}} \right|_{\mathcal{F}_t} = e^{cX_t-\psi(c)t}.
\end{equation}
Under $\mathbb{P}^c$, the process $X$ is still a spectrally negative L\'evy process and we mark its Laplace exponent and scale functions with the subscript $c$. From $\psi_c(\lambda)=\psi(\lambda+c)-\psi(c)$ for $\lambda \geq 0$ we get by taking Laplace transforms

\[ W^{(q)}_c(x)=e^{-cx} W^{(q+\psi(c))}(x) \]
for all $q \geq 0$ and, similarly
\[Z^{(q)}_c(x)=1+q\int_0^x W^{(q)}_c (y)\, \mbox{d}y\]
\subsection{Reviewing the McKean stochastic game}\label{subsec_review}
First consider the McKean optimal stopping problem (or American put option) with value function $U$, i.e.
\[U(x)=\sup_\tau\mathbb{E}_x[e^{-q\tau}(K-e^{X_\tau})^+].\]
We recall the solution to this problem as it appears in \cite{Chan} (see also \cite{Mordecki}):

\begin{theorem}\label{avram}
For the McKean optimal stopping problem under (\ref{ass}) we have
\[
U(x) = K Z^{(q)}(x-k^*) - e^x Z_1^{(q-\psi(1))}(x-k^*),
\]
where
\[
e^{k^*} = K\frac{q}{\Phi(q)}\frac{\Phi(q)-1}{q-\psi(1)},
\]
which is to be understood in the limiting sense when $q=\psi(1)$, in other words, $e^{k^*} = K \psi(1)/\psi'(1)$. An optimal stopping time is given by $\tau^*=\inf\{t>0 : X_t < k^*\}$.
\end{theorem}

Next we recall the main result from BK on a saddle point and the value function for the McKean game:
\begin{theorem}\label{mainthrm}Consider the McKean stochastic game under the assumption (\ref{ass}) and recall that $\bar{\delta}=U(\log K)$.
\begin{itemize}
\item[(i)] If $\delta \geq \bar{\delta}$, then a stochastic saddle point is given by $\tau^*$ from Theorem \ref{avram} and $\sigma^*=\infty$, in which case $V=U.$
\item[(ii)] If $\delta< \bar{\delta}$, a stochastic saddle point is given by the pair
\[
\tau^*= \inf\{t>0 : X_t <x^*\} \text{ and  }\sigma^*=\inf\{t> 0 : X_t \in [\log K, y^*]\},
\]
where $x^*$ uniquely solves
\begin{equation}
Z^{(q)}(\log K -x) -  Z_1^{(q-\psi(1))} (\log K - x)= \frac{\delta}{K},
\label{howtofindx*}
\end{equation}
$x^*> k^*$ (the optimal level of the corresponding McKean optimal stopping problem in Theorem \ref{avram}) and $y^* \geq \log K$.

Furthermore,
\begin{equation}V(x) = KZ^{(q)}(x - x^*) - e^xZ_1^{(q-\psi(1))}(x- x^*)\label{equationV}\end{equation}
for $x\leq \log K$ and if $y^*=\log K$ then for any $x \in \mathbb{R}$
\[ V(x) = KZ^{(q)}(x-x^*)-e^x Z_1^{(q-\psi(1))}(x-x^*) + \alpha e^{\Phi (q)(\log K-x^*)} W^{(q)}(x-\log K), \]
where
\[ \alpha = e^{x^*} \frac{q-\psi(1)}{\Phi (q)-1}-\frac{qK}{\Phi (q)}, \]
which is to be understood in the limiting sense when $q=\psi(1)$, i.e. $\alpha=e^{x^*} \psi'(1)-K\psi(1)$.
\end{itemize}
\end{theorem}
Hence a saddle point exists, and consists of the first hitting time of $(-\infty,x^*]$ for the maximiser and of the first hitting time of $[\log K,y^*]$ for the minimiser. Note that when $\delta<\bar{\delta}$ we know that the value of the McKean is not everywhere equal to that of its one-player counterpart as $V(\log K)=\delta<U(\log K)$ in that case.
Furthermore, equation (\ref{howtofindx*}) gives us a characterisation of $x^*$, but we know only little about $y^*$.

In BK, the issue whether $y^*=\log K$ or $y^*>\log K$ was only answered when $X$ has no Gaussian component:

\begin{theorem}
\label{mainthrmII}
Suppose in Theorem \ref{mainthrm} that $\delta<  \bar{\delta}$.
If $X$ has no Gaussian component, then $y^*>\log K$ and necessarily $\Pi(-\infty, \log K -y^*)>0$.
\end{theorem}

\begin{remark}\label{rem_Kees1}

We now discuss the intuitive interpretation of these results.

Due to the fact that $X$ has no positive jumps, the choice of $y^*$ has no influence when $X$ is started at a value to the left of $\log K$ as the minimiser would stop as soon as the process hits $\log K$.

However, the situation is different when starting from any $X_0>\log K$, the minimiser could either stop right away and pay $\delta$ to the maximiser, or wait a short time. From the minimiser's point of view, the latter decision has the advantage of profiting from the discounting, but the disadvantage of the risk that a (large) negative jump could bring $X$ (far) below $\log K$, where a higher payoff than (discounted) $\delta$ can be claimed by the maximiser. The closer $X_0$ is chosen to $\log K$, the more dominant the disadvantage becomes and hence the exercise region for the minimiser can take the form of an interval $[\log K,y^*]$.

When $X$ is a Brownian motion it is obvious that we have $y^*=\log K$ for any $\delta \in (0,\bar{\delta}]$ (see also \cite{Kyprianou04}) as then the process can only get below $\log K$ by hitting it first.

The above Theorem \ref{mainthrmII} tells us that the other extreme case, namely $y^*>\log K$ for any $\delta \in (0,\bar{\delta}]$, i.e. the disadvantage of waiting being dominant for the minimiser, occurs whenever $X$ has no Gaussian component.

The interesting question is what happens when $X$ has a Gaussian component \emph{and} negative jumps as then there is a trade-off between the possibility for the process to enter the region $(-\infty,\log K]$ continuously (which can only occur when $\sigma_X>0$), leading to a relatively small pay-off and the possibility of passing $\log K$ by a jump potentially leading to a larger pay-off.
 It turns out that for $\delta$ large enough, when stopping immediately is relatively expensive, the Gaussian part `wins' in the sense that $y^*=\log K$ (the minimiser is happy to take the risk of the process jumping to a less-favourable region), while for $\delta$ small enough, when stopping immediately has become cheaper, the negative jumps `win' in the sense that $y^*>\log K$, see Theorem \ref{thm_Erik_main} below. In fact, in Figure 4 at the end of this paper we plot the negative relationship between $y^*$ and $\sigma_X$ in the jump-diffusion case.
\end{remark}

\section{Single point or interval when $X$ has a Gaussian part $\sigma_{X}>0$}\label{sec_main_res}

Throughout this section we assume that condition (\ref{ass}) holds. Recall that $T_K := \inf \{ t>0 \, | \, X_t = \log K \}$. 
Consider the following function
\begin{equation}f_{\delta}(x)=\sup_\tau \mathbb{E}_x[e^{-q \tau}(K-e^{X_\tau})\mathbf{1}_{\{\tau\leq T_K\}}+\delta e^{-q T_K}\mathbf{1}_{\{T_K<\tau\}}], \label{auxosp}
\end{equation}
i.e. the optimal value for the maximiser provided the minimiser only exercises when $X$ hits $\log K$.

We first prove the following technical result.
\begin{lemma}\label{hulp}
Suppose $\sigma_{X}>0$ and $0<\delta\leq \bar{\delta}=U(\log K)$.
The function $f_{\delta}$ is differentiable on $\mathbb{R}\backslash \{\log K\}$.
Furthermore, $f_{\delta}=V$ on $(-\infty,\log K],$ $f_{\delta}\geq V$ on $\mathbb{R}$ and $f_{\delta}'(\log K+)$ is a strictly decreasing continuous function of $\delta$.
\end{lemma}
\begin{proof}
Let $\delta \in (0,\bar{\delta}]$.
Due to Theorem \ref{mainthrm} and the absence of positive jumps we have for $x\leq \log K$
\begin{eqnarray*}
V(x)&=&\mathbb{E}_x[e^{-q\tau_{x^*(\delta)}^-}(K-e^{X_{\tau_{x^*(\delta)}^-}})\mathbf{1}_{\{\tau_{x^*(\delta)}^-<T_K\}}+\delta e^{-qT_K}\mathbf{1}_{\{T_K<\tau_{x^*(\delta)}^-\}}]\\
&=&\sup_{\tau}\mathbb{E}_x[e^{-q\tau}(K-e^{X_{\tau}})\mathbf{1}_{\{\tau<T_K\}}+\delta e^{-qT_K}\mathbf{1}_{\{T_K<\tau\}}]\\
&=&f_{\delta}(x).
\end{eqnarray*}
Also, for any $x\in\mathbb{R}$
\begin{eqnarray*}f_{\delta}(x)&=&\sup_\tau \mathbb{E}_x[e^{-q \tau}(K-e^{X_\tau})\mathbf{1}_{\{\tau\leq T_K\}}+\delta e^{-q T_K}\mathbf{1}_{\{T_K<\tau\}}]\\
&\geq &\inf_\sigma\sup_\tau \mathbb{E}_x[e^{-q \tau}(K-e^{X_\tau})\mathbf{1}_{\{\tau\leq \sigma\}}+\delta e^{-q \sigma}\mathbf{1}_{\{\sigma<\tau\}}]\\
&=&V(x).
\end{eqnarray*}
In fact, since stopping is not optimal on $(\log K,\infty)$ as the lower pay-off function is zero there, we deduce that we have for all $x\in\mathbb{R}$
\begin{equation}\label{Kees_28jul1}
f_{\delta}(x)= \mathbb{E}_x[e^{-q \tau_{x^*(\delta)}^-}(K-e^{X_{\tau_{x^*(\delta)}^-}})\mathbf{1}_{\{\tau_{x^*(\delta)}^-\leq T_K\}}+\delta e^{-q T_K}\mathbf{1}_{\{T_K<\tau_{x^*(\delta)}^-\}}].
\end{equation}
Now, let $\delta_2>\delta_1>c$ for some $c>0$. From the defintion of $f_{\delta}$ in (\ref{auxosp}) we find
\begin{eqnarray*}f_{\delta_2}(x)-f_{\delta_1}(x)&=&\sup_\tau \mathbb{E}_x[e^{-q \tau}(K-e^{X_\tau})\mathbf{1}_{\{\tau\leq T_K\}}+\delta_2 e^{-q T_K}\mathbf{1}_{\{T_K<\tau\}}]\\
&&-\sup_\tau \mathbb{E}_x[e^{-q \tau}(K-e^{X_\tau})\mathbf{1}_{\{\tau\leq T_K\}}+\delta_1 e^{-q T_K}\mathbf{1}_{\{T_K<\tau\}}]\\
&\leq&(\delta_2-\delta_1)\sup_\tau \mathbb{E}_x[e^{-q T_K}\mathbf{1}_{\{T_K<\tau\}}]\\
&\leq & (\delta_2-\delta_1) \mathbb{E}_x[e^{-q \tau_{\log K}^-}],
\end{eqnarray*}
from which it follows that (the equality by (\ref{K_30sept1}))
\begin{eqnarray*}\frac{f_{\delta_2}(\log K+\varepsilon)-\delta_2}{\varepsilon}-\frac{f_{\delta_1}(\log K+\varepsilon)-\delta_1}{\varepsilon}&\leq& (\delta_2-\delta_1) \frac{\mathbb{E}_{\log K+\varepsilon}[e^{-q \tau_{\log K}^-}]-1}{\varepsilon}\\
&=&(\delta_2-\delta_1) \left( \frac{Z^{(q)}(\varepsilon)-1}{\varepsilon} - \frac{q}{\Phi(q)} \frac{W^{(q)}(\varepsilon)}{\varepsilon} \right).
\end{eqnarray*}

Since $f_{\delta}$ is a differentiable function on $[\log K,\infty)$ (see equation (27) in BK together with (\ref{Kees_28jul1})) and using $Z^{(q)\prime}(0)=W^{(q)}(0)=0$, $W^{(q)\prime}(0+)=2/\sigma_{X}^2$ we deduce that
\begin{equation}\label{K_30sept3}
f_{\delta_2}'(\log K+)-f_{\delta_1}'(\log K+)\leq -\frac{2q}{\sigma_{X}^2\Phi(q)}(\delta_2-\delta_1),
\end{equation}
showing that $f_{\delta}'(\log K+)$ is strictly decreasing in $\delta$.

Also,  using (\ref{auxosp}) and the fact that $\tau_{x^*(\delta_1)}^-$ is a feasible strategy also when $\delta=\delta_2$, it holds that
\begin{eqnarray*}f_{\delta_2}(x)-f_{\delta_1}(x)&\geq &\mathbb{E}_x[e^{-q \tau_{x^*(\delta_1)}^-}(K-e^{X_{\tau_{x^*(\delta_1)}^-}})\mathbf{1}_{\{\tau_{x^*(\delta_1)}^-\leq T_K\}}+\delta_2 e^{-q T_K}\mathbf{1}_{\{T_K<\tau_{x^*(\delta_1)}^-\}}]\\
&&-\mathbb{E}_x[e^{-q \tau_{x^*(\delta_1)}^-}(K-e^{X_{\tau_{x^*(\delta_1)}^-}})\mathbf{1}_{\{\tau_{x^*(\delta_1)}^-\leq T_K\}}+\delta_1 e^{-q T_K}\mathbf{1}_{\{T_K<\tau_{x^*(\delta_1)}^-\}}]\\
&=&(\delta_2-\delta_1)\mathbb{E}_x[e^{-q T_K}\mathbf{1}_{\{T_K<\tau_{x^*(\delta_1)}^-\}}]\\
&\geq&(\delta_2-\delta_1)\mathbb{E}_x[e^{-q T_K}\mathbf{1}_{\{T_K<\tau_{x^*(c)}^-\}}],
\end{eqnarray*}
where the final inequality follows from the observation that $x^*(\delta)$ is decreasing in $\delta$ and that $\delta_1>c$.
Note that $x^*(c)<\log(K-c)$ since $V(x)$ is strictly decreasing in $x\in(-\infty,\log K]$ for any $\delta>0$ and thus
\begin{eqnarray*}&&\hspace{-3cm}\frac{f_{\delta_2}(\log K+\varepsilon)-\delta_2}{\varepsilon}-\frac{f_{\delta_1}(\log K+\varepsilon)-\delta_1}{\varepsilon}\\
&\geq& (\delta_2-\delta_1) \frac{\mathbb{E}_{\log K+\varepsilon}[e^{-q T_K}\mathbf{1}_{\{T_K<\tau_{x^*(c)}^-\}}]-1}{\varepsilon}\\
&=&(\delta_2-\delta_1)\frac{W^{(q)}(\log K+\varepsilon-x^*(c))-W^{(q)}(\log K -x^*(c))}{\varepsilon W^{(q)}(\log K -x^*(c))}\\
&&-(\delta_2-\delta_1)e^{\Phi(q)(\log K -x^*(c))}\frac{W^{(q)}(\varepsilon)}{\varepsilon W^{(q)}(\log K-x^*(c))}.
\end{eqnarray*}
because of Lemma 12 in BK.
It follows that
\[f_{\delta_1}'(\log K+)-f_{\delta_1}'(\log K+)\geq (\delta_2-\delta_1)\frac{\sigma_{X}^2W^{(q)\prime}(\log K-x^*(c))- 2e^{\Phi(q)(\log K -x^*(c))}}{\sigma_{X}^2W^{(q)}(\log K-x^*(c))}.\]
Since $c$ is arbitrary, we conclude from this inequality together with (\ref{K_30sept3}) that $f_{\delta}'(\log K+)$ is indeed continuous in $\delta$ for any $\delta>0$.
\end{proof}

Now we are ready to prove our main result, extending Theorem \ref{mainthrm}:

\begin{theorem}\label{thm_Erik_main}
Suppose $\sigma_{X}>0$. When $\Pi\neq 0$,  then there exists a unique $\delta_0 \in (0,\bar{\delta})$ such that an optimal stopping time for the minimiser is given by
$T_K$ (i.e. $y^*(\delta)=\log K$) when $\delta \in [\delta_0,\bar{\delta}]$ and by $T_{[\log K, y^*(\delta)]}$ for some $y^*(\delta)>\log K$ when $\delta \in (0,\delta_0)$.
\end{theorem}

\begin{proof}
Let $\sigma_{X}>0$ and suppose $\Pi\neq 0$. We know from Theorem \ref{mainthrm} that the stopping region for the minimiser is of the form $[\log K,y^*]$ for some $y^*\geq \log K$. We claim that setting $\delta_0$ equal to the unique zero of $f_{\delta}'(\log K+)$ on $(0,\bar{\delta})$ yields the result.

First let us show that this unique zero indeed exists. For $\delta=\bar{\delta}$ it holds that $f_{\delta}'(\log K+)=U'(\log K)<0$ (cf. Theorem \ref{avram}). Using Lemma \ref{hulp}, it suffices to show that there exists some $\delta>0$ such that $f_{\delta}'(\log K+)>0$. We argue by contradiction, so, again using Lemma \ref{hulp}, suppose that $f_{\delta}'(\log K+)<0$ for all $\delta>0$. This implies that for each $\delta>0$ there exists some $\varepsilon>0$ such that $f_{\delta}(x)<f_{\delta}(\log K)=\delta$ for all $x\in(\log K,\log K+\varepsilon].$ Since $V \leq f_{\delta}$ (Lemma \ref{hulp}) we deduce that $V(x)<\delta=(K-e^x)^++\delta$ for all $x\in(\log K, \log K +\varepsilon)$, hence $y^*=\log K$ and in fact $V=f_{\delta}$ (by (\ref{auxosp})).

But plugging $\tau_{\log K/2}^-$ in the rhs of (\ref{auxosp}) yields
\[f_{\delta}(x)\geq K/2\mathbb{E}_x[e^{-q\tau_{\log K/2}^-}\mathbf{1}_{\{\tau_{\log K/2}^-<T_K\}}].\]
This lower bound is strictly positive for $x>\log K$ since $\Pi\neq 0$ and does not depend on $\delta$. Hence for $\delta$ small enough we deduce the existence of some $x>\log K$ such that $f_{\delta}(x)>\delta$, which contradicts with $f_{\delta}(x)=V(x) \leq \delta$ on $[\log K,\infty)$.

Next for the optimal stopping time of the minimiser. For $\delta>\delta_0$ the same reasoning as above yields $y^*=\log K$.

For the case $\delta=\delta_0$ we note that for any fixed $x$ the function $f_{\delta}(x)$ is continuous in $\delta$, as is easily seen from (\ref{auxosp}). Hence
\[f_{\delta_0}(x)=\lim_{\delta\downarrow \delta_0}f_{\delta}(x)\leq (K-e^x)^++\delta_0,\]
from which we can deduce that we still have $y^*=\log K$. Finally, let $\delta<\delta_0$. Again much as above, we then have that $f_{\delta}'(\log K+)>0$ and thus there exist $x>\log K$ for which $f_{\delta}(x)>\delta=(K-e^x)^++\delta$. Since trivially $V$ is bounded above by this upper payoff function, it cannot be true that $f_{\delta}=V$ and thus it can also not be true that $y^*=\log K$, so we indeed arrive at $y^*>\log K$.

\end{proof}

\begin{remark} From the proof of the above Theorem \ref{thm_Erik_main} we see that this result is essentially due to the upper payoff function $(K-e^x)^++\delta$ having a kink at the point where it first touches the value function as $\delta$ decreases (namely $\log K$). That is, if we would only slightly alter the upper payoff function on an environment of $\log K$ so it would have a continuous derivative, we should expect the optimal stopping time for the minimiser to be $T_{[y_1^*,y_2^*]}$ with $y_1^*<\log K<y_2^*$ for \emph{all} $\delta \in (0,\bar{\delta})$ and \emph{any} spectrally negative L\'evy process $X$.
\end{remark}

Next we provide expressions that complement those from Theorem \ref{mainthrm}. Recall that equation(\ref{equationV}) in Theorem \ref{mainthrm} provides us with a formula for $V$ on $(-\infty,\log K]$, so we can make use of the following function:

\begin{equation}\label{wdelta}w_\delta(x)=\left\{\begin{array}{ll}
V(x)&\mbox{for $x<\log K$}\\
\delta&\mbox{for $x\geq \log K$}.
\end{array}
\right.
\end{equation}

\begin{theorem}\label{thm_Erik2}
Suppose $\Pi \not= 0$. We have the following.
\begin{itemize}
\item[(i)] Suppose $\sigma_{X}>0$. Then $\delta_0$ is the unique solution on $(0,\bar{\delta})$ to the equation in $\delta$:
\[ \int_{t<0}\int_{u<t}(w_\delta(t+\log K)-\delta)e^{-\Phi(q)(t-u)}\Pi(du)dt=\frac{\delta q}{\Phi(q)}. \]
\item[(ii)] Suppose $y^*>\log K$ (i.e. $\sigma_{X}>0$ and $\delta <\delta_0$, or $\sigma_{X}=0$ and $\delta<\bar{\delta}$). Then $y^*$ is the unique solution on $(\log K,\infty)$ to the equation in $y$:
\begin{equation}\label{Kees_29jul2}
\int_{t<0}\int_{u<t}(w_\delta(t+y)-\delta)e^{-\Phi(q)(t-u)}\Pi(du)dt=\frac{\delta q}{\Phi(q)}.
\end{equation}
Furthermore, $V(x)=\delta$ for $x \in [\log K,y^*]$ and for $x \in (y^*,\infty)$:
\begin{equation}\label{29jul3}
V(x) = \delta Z^{(q)}(x-y^*) - \int_{t<0}\int_{u<t}(w_\delta(t+y^*)-\delta) W^{(q)}(x-y^*-t+u) \Pi(du)dt.
\end{equation}
\end{itemize}
\end{theorem}

\begin{proof} First we introduce the function
\begin{equation}h(x,y):=\mathbb{E}_x[e^{-q\tau_{y}^-}w_\delta(X_{\tau_y^-})] \label{h}\end{equation}
for $x>y \geq \log K$. Observe that by the lack of positive jumps, $h(.,y)$ is the optimal value the maximiser can obtain when the minimiser chooses as stopping region $[\log K,y]$. Hence in particular $V(x)=h(x,y^*)$.

Denote by $u^{(q)}(s,t)$ the resolvent density of $X$ started at $s>0$ and killed at first passage below $0$. Invoking the compensation formula (see e.g. Theorem 4.4 in \cite{Kyprianou06}) leads to
\begin{eqnarray*}
h(x,y)&=&\delta \mathbb{E}_x[e^{-q\tau_y^-}]+\mathbb{E}_x[e^{-q\tau_{y}^-}(w_\delta(X_{\tau_y^-})-\delta)\mathbf{1}_{\{X_{\tau_y^-}<\log K\}}]\\
&=&\delta \mathbb{E}_x[e^{-q\tau_y^-}]+\int_{t<\log K-y}\int_{u<t}(w_\delta(t+y)-\delta)u^{(q)}(x-y,t-u)\Pi(du)dt\\
&=&\delta \mathbb{E}_x[e^{-q\tau_y^-}]+\int_{t<0}\int_{u<t}(w_\delta(t+y)-\delta)u^{(q)}(x-y,t-u)\Pi(du)dt,
\end{eqnarray*}
where the final equality is due to the fact that $w_\delta=\delta$ on $[\log K,y]$.
We know that  (see e.g. Theorem 8.1 and Corollary 8.8 in \cite{Kyprianou06} respectively)
\[\mathbb{E}_x[e^{-q\tau_y^-}]=Z^{(q)}(x-y)-\frac{q}{\Phi(q)}W^{(q)}(x-y)\]
and
\[u^{(q)}(s,t)=e^{-\Phi(q)t}W^{(q)}(s)-W^{(q)}(s-t),\]
hence
\begin{eqnarray}
h(x,y)&=&
\int_{t<0}\int_{u<t}(w_\delta(t+y)-\delta)(e^{-\Phi(q)(t-u)}W^{(q)}(x-y)-W^{(q)}(x-y-t+u))\Pi(du)dt\nonumber
\\
&&+\delta(Z^{(q)}(x-y)-\frac{q}{\Phi(q)}W^{(q)}(x-y)).\label{compensation2}
\end{eqnarray}

Furthermore, when $X$ is of unbounded variation we can compute for $x>y$
\begin{eqnarray*}
\frac{\partial}{\partial x}h(x,y)&=&\delta(qW^{(q)}(x-y)-\frac{q}{\Phi(q)}W^{(q)^\prime}(x-y))\\
&&\hspace{-2cm}+\int_{t<0}\int_{u<t}(w_\delta(t+y)-\delta)(e^{-\Phi(q)(t-u)}W^{(q)^\prime}(x-y)-W^{(q)^\prime}(x-y-t+u))\Pi(du)dt.
\end{eqnarray*}
and we can let $x \downarrow y$ to arrive at
\begin{equation}\label{Kees_29jul1}
\frac{\partial}{\partial x}h(y+,y) = \left(\int_{t<0}\int_{u<t}(w_\delta(t+y)-\delta)e^{-\Phi(q)(t-u)}\Pi(du)dt-\frac{q\delta}{\Phi(q)}\right)W^{(q)^\prime}(0+).
\end{equation}

Ad (i). Recall the function $f_{\delta}$ as defined in (\ref{auxosp}), and recall in particular from the proof of Lemma \ref{hulp} that $\delta_0$ is the unique $\delta \in (0,\bar{\delta})$ for which $f_{\delta}'(\log K+)=0$. Furthermore, note that $f_{\delta}(x)=h(x,\log K)$ for $x>\log K$, since both sides equal the optimal value the maximiser can obtain when the minimiser only stops when $X$ hits $\log K$. Combining these observations with (\ref{Kees_29jul1}) and $W^{(q)^\prime}(0+)=2/\sigma_{X}^2 \not= 0$ yields the result.

Ad (ii). We first consider the case when $X$ is of bounded variation. We know from Theorem 4 in BK that we have continuous fit, i.e. $V(y^*+)=\delta$. Since the integrand in (\ref{compensation2}) is bounded and equal to zero for $t<\log K-y$ we can take the limit inside the integrals to deduce that
\[h(y+,y)=\delta-\frac{q\delta}{\mathrm{d}\Phi(q)} +\frac{1}{\mathrm{d}}\int_{t<0}\int_{u<t}(w_\delta(t+y)-\delta)e^{-\Phi(q)(t-u)}\Pi(du)dt,\]
so using $V(y^*+)=h(y^*+,y^*)$ it follows that $y^*$ indeed solves (\ref{Kees_29jul2}). For uniqueness, the function $w_\delta=V$ is strictly decreasing on $(-\infty,\log K]$ and $\delta=V(y^*)=h(y^*+,y^*).$ Since $q>0$, the minimiser would not stop at points in $[\log K,\infty]$ from which the process cannot jump into $(-\infty,\log K)$ and thus $\log K-y^*>l:=\sup\{x:\Pi(-\infty,x)=0\}.$ Combining these observations imply that $h(y+,y)$ is a strictly decreasing function on $[\log K,\log K-l]$.

Next consider the case that $X$ is of unbounded variation. Now Theorem 4 in BK tells us that we have smooth fit, i.e. $V'(y^*+)=0$. Using $V(x)=h(x,y^*)$ together with (\ref{Kees_29jul1}) yields again that $y^*$ solves (\ref{Kees_29jul2}), uniqueness follows in the same way as in the previous paragraph.

Finally, (\ref{29jul3}) is readily seen from $V(x)=h(x,y^*)$, (\ref{compensation2}) and the fact that $y^*$ satisfies (\ref{Kees_29jul2}).

\end{proof}

We conclude this section with some properties of $y^*$ as a function of $\delta$. Note that by spectral negativity, $\Pi \not=0$ implies $\sup \{ x \, : \, \Pi(-\infty,x)=0 \} < 0$.

\begin{theorem} Suppose $\Pi \not=0$. Then $y^*(\delta)$ is continuous and decreasing as a function of $\delta$, with $y^*(\bar{\delta}-)=\log K$ if $\sigma_{X}=0$ (resp. $y^*(\delta_0-)=\log K$ if $\sigma_{X}>0$) and $y^*(0+)=\log K-\sup \{ x \, : \, \Pi(-\infty,x)=0 \}$.
\end{theorem}

\begin{proof} We write $V_{\delta}$ to stress the dependence of the value function on $\delta$. Continuity of $y^*(\delta)$ is clear as the above Theorem \ref{thm_Erik2} (ii) and the fact that $w_{\delta}$ is continuous in $\delta$ (see the argument for continuity of $\delta \mapsto V_{\delta}$ below) allow to apply the implicit function theorem.

To see that it is decreasing it suffices to show that $\delta \mapsto V_{\delta}(x)-\delta$ is decreasing. For this, take $\delta_1<\delta_2$ and let $(\tau^*_1,\sigma^*_1)$ denote the saddle point when $\delta=\delta_1$. Then $V_{\delta_1}$ is the value when the supremum over all $\tau$ is taken in the expected pay-off corresponding to the pair $(\tau,\sigma^*_1)$. As $\sigma^*_1$ is also feasible for the minimiser when $\delta=\delta_2$ we have that $V_{\delta_2}$ is bounded above by the value when the supremum over all $\tau$ is taken in the expected pay-off corresponding to the pair $(\tau,\sigma^*_1)$. This yields

\begin{eqnarray}
V_{\delta_2}(x)-V_{\delta_1}(x) & \leq & \sup_{\tau} \mathbb{E}_x [ e^{-q \sigma^*_1}((K-e^{X_{\sigma^*_1}})^+ +\delta_2) \mathbf{1}_{\{ \sigma^*_1<\tau \}} \nonumber\\
 & & \quad - e^{-q \sigma^*_1} ((K-e^{X_{\sigma^*_1}})^+ +\delta_1) \mathbf{1}_{\{ \sigma^*_1<\tau \}} ] \nonumber\\
 & \leq & \delta_2-\delta_1, \label{Kees_30jul5}
\end{eqnarray}
as required.

Next, by the monotonicity the limits mentioned in the theorem exist. First we show $y^*(0+)=\log K-l$, where $l:=\sup \{ x \, : \, \Pi(-\infty,x)=0 \}$. Suppose we had $y^*(0+)<\log K-l$, then for some $x_1 \in (y^*(0+),\log K-l)$ and any $\delta>0$ we have $\mathbb{P}_{x_1}(\tau^-_{\log K/2}<T_{[\log K,y^*(\delta)]}) \geq \mathbb{P}_{x_1}(\tau^-_{\log K/2}<T_{[\log K,y^*(0+)]})>0$. So, starting from $x_1$, if the maximiser chooses $\tau^-_{\log K/2}$ he ensures a strictly positive value, independent of $\delta$. But this of course contradicts with $V_{\delta}(x_1) \leq \delta \downarrow 0$ as $\delta \downarrow 0$. If we had $y^*(0+)>\log K-l$, then for some $x_2 \in (\log K-l,y^*(0+))$ we have for $\delta$ small enough $x_2 \leq y^*(\delta)$ and consequently $V_{\delta}(x_2)=\delta$. But the minimiser can do better, that is in fact we have $V_{\delta}(x_2)<\delta$, as is easily seen. Namely, the minimiser can choose $T_{[\log K,\log K-l]}$, so that starting from $x_2>\log K-l$ the maximiser can at most get discounted $\delta$, the discount factor being strictly less than $1$ since $q>0$ and $X$ is right continuous.

Next suppose $\sigma_{X}>0$ and let us show that $y^*(\delta_0-)=\log K$. Suppose we had $y^*(\delta_0-)>\log K$. Note that for any $x$, $\delta \mapsto V_{\delta}(x)$ is continuous, since for $\delta_1<\delta_2$ trivially $V_{\delta_2}(x) \geq V_{\delta_1}(x)$ and (\ref{Kees_30jul5}). So for $\log K<x_1<x_2<y^*(\delta_0-)$ it would follow that $V_{\delta}(x_1)-V_{\delta}(x_2) \to V_{\delta_0}(x_1)-V_{\delta_0}(x_2)=\delta_0-\delta_0=0$ as $\delta \downarrow \delta_0$. But the difference $V_{\delta}(x_1)-V_{\delta}(x_2)$ does not vanish as $\delta \downarrow \delta_0$, as follows easily from the homogeneity of $X$. More precisely, denoting by $(\tau^*_1,\sigma^*_1)$ resp. $(\tau^*_2,\sigma^*_2)$ the saddle point when starting from $x_1$ resp. $x_2$, similar arguments as the ones leading to (\ref{Kees_30jul5}) yield in this case

\[ V_{\delta}(x_1) \geq \mathbb{E} [ e^{-q \tau^*_2} (K-e^{x_1 + X_{\tau^*_2}})^+ \mathbf{1}_{\{ \tau^*_2 \leq \sigma^*_1 \}} +  e^{-q \sigma^*_1} ((K-e^{x_1 + X_{\sigma^*_1}})^+  +\delta)\mathbf{1}_{\{  \sigma^*_1 <\tau^*_2 \}} ] \]
and

\[ V_{\delta}(x_2) \leq \mathbb{E} [ e^{-q \tau^*_2} (K-e^{x_2 + X_{\tau^*_2}})^+ \mathbf{1}_{\{ \tau^*_2 \leq \sigma^*_1 \}} +  e^{-q \sigma^*_1} ((K-e^{x_2 + X_{\sigma^*_1}})^+  +\delta)\mathbf{1}_{\{  \sigma^*_1 <\tau^*_2 \}} ], \]
thus

\begin{equation}\label{Kees_31jul1}
V_{\delta}(x_1) - V_{\delta}(x_2) \geq \mathbb{E} [ e^{-q \kappa} ( (K-e^{x_1 + X_{\kappa}})^+ - (K-e^{x_2 + X_{\kappa}})^+) ]
\end{equation}
where $\kappa=\sigma^*_1 \wedge \tau^*_2=\inf \{ t>0 \, | \, X_t = \log K-x_1 \} \wedge \inf \{ t > 0 \, | \, X_t < x^*(\delta)-x_2 \}$. Clearly, since $x^*(\delta) \leq \log K$ and $x_1<x_2$ the rhs of (\ref{Kees_31jul1}) is strictly positive iff $\mathbb{P}(\tau^*_2<\sigma^*_1)>0$. Even after taking the limit for $\delta \downarrow \delta_0$ this probability is positive on account of $\Pi \not=0$.

Finally, $y^*(\bar{\delta}-)=\log K$ when $\sigma_{X}=0$ can be shown by the same arguments, taking into account here one has $\sigma^*=\infty$ for $\delta>\bar{\delta}$.
\end{proof}

\section{Jump-diffusion case}\label{sec_jump_diff}

In this section we translate the general results from the previous Section \ref{sec_main_res} to the particular case of a jump-diffusion with downwards directed, exponentially distributed jumps. In this case, which is quite popular in practical applications in finance e.g. due to its tractable nature, the expressions become much more explicit. In particular a formula exists that expresses $y^*$ explicit in terms of $x^*$, cf. Proposition \ref{K_1mei_3} (iv).

For the sequel we set

\begin{equation}\label{K_6apr_1}
X_t = \sigma_{X} W_t + \mu t - \sum_{i=1}^{N_t} \xi_i, \quad t \geq 0,
\end{equation}
where $\sigma_{X}>0$, $\mu \in \mathbb{R}$, $N$ is a Poisson process with intensity $\lambda>0$ counting the jumps and $(\xi_i)_{i \geq 0}$ is an iid sequence of random variables following an exponential distribution with parameter $\theta>0$.

The following Proposition \ref{K_7apr_1} states formulas for the scale functions in this jump-diffusion case (recall $\mathbb{P}^c$ as defined in (\ref{K_18sept1})):

\begin{proposition}\label{K_7apr_1}
Let $c,r \geq 0$. We have the following for $X$ given by (\ref{K_6apr_1}) under $\mathbb{P}^c$.
\begin{itemize}
\item[(i)] The Laplacian is given by

\[ \psi_c(z) = \psi(z+c) - \psi(c) = \frac{\sigma_{X}^2}{2} z^2 +(\sigma_{X}^2 c+ \mu) z - \frac{\lambda \theta z}{(\theta+z+c)(\theta+c) }. \]
The function $z \mapsto \psi_c(z)-r$ has three zeros $\beta_1(c,r)<-\theta-c<\beta_2(c,r) \leq \beta_3(c,r)$, with $\beta_2(c,r)<0<\beta_3(c,r)$ if $r>0$; $\beta_2(c,r)=0<\beta_3(c,r)$ if $r=0$ and $\psi_c'(0) \leq 0$; $\beta_2(c,r)<0=\beta_3(c,r)$ if $r=0$ and $\psi_c'(0) \geq 0$.

\item[(ii)] In particular, if $r=\psi(1)>0$ we have

\begin{multline*}
\beta_{1,2}(0,r) = - \left( \frac{\theta}{2} + \frac{r}{\sigma_{X}^2} + \frac{\lambda}{\sigma_{X}^2 (\theta+1)} \right) \pm \sqrt{\left( \frac{\theta}{2} + \frac{r}{\sigma_{X}^2} + \frac{\lambda}{\sigma_{X}^2 (\theta+1)} \right)^2 - \frac{2r\theta}{\sigma_{X}^2}} \\
\mbox{and} \quad \beta_{3}(0,r)=1.
\end{multline*}
\end{itemize}
Define for $i=1,2,3$ the constants

\[ C_i(c,r) = \frac{2 (\theta + c + \beta_i(c,r))}{\sigma_{X}^2 \prod_{j \not= i} ( \beta_j(c,r)-\beta_i(c,r) )}. \]
We have the following formulas for the scale functions $W^{(r)}_c$ and $Z^{(r)}_c$ on $[0,\infty)$.

\begin{itemize}
\item[(iii)] If $\beta_2(c,r)\not=0$ or $\beta_3(c,r)\not=0$ we have

\[ W^{(r)}_c(x) = \sum_{i=1}^3 C_i(c,r) e^{\beta_i(c,r) x}, \]
otherwise (necessarily $r=0$) we have

\[ W^{(0)}_c(x) = \frac{2}{\sigma_{X}^2 \beta_1(c,0)} \left( (1-c-\theta) e^{\beta_1(c,0) x} -(\theta+c) x + \theta+c-1 \right). \]

\item[(iv)] If $r>0$ we have

\[ Z^{(r)}_c(x) = r \sum_{i=1}^3 \frac{C_i(c,r)}{\beta_i(c,r)} e^{\beta_i(c,r) x},\]
while $Z^{(0)}_c(x)=1$.
\end{itemize}
\end{proposition}

\begin{proof} Follows from the definitions (\ref{K_18sept2}) and (\ref{K_18sept3}) by some elementary calculations.  Also, see e.g. \cite{Avram04}.
\end{proof}

In Propositions \ref{K_22apr_1} and \ref{K_1mei_3} below we assume for simplicity that $q>0$ and $q=\psi(1)$, i.e. we set $\mu :=q-\sigma_{X}^2/2+\lambda/(\theta+1)$. (Note that condition (\ref{ass}) is met). This means that $\mathbb{P}$ is a so-called risk neutral measure in the sense that the discounted price process $(e^{X_t-qt})_{t \geq 0}$ is a $\mathbb{P}$-martingale, as required in a financial modelling context. (However the reader should have no difficulties translating the upcoming formulas to the situation for any $q \in [0,\psi(1)]$ if required.) Note that the above Proposition \ref{K_7apr_1} (ii) gives explicit formulas for the roots $\beta_i(0,q)$ in this case.

First we turn to formulas for the McKean optimal stopping problem (cf. Theorem \ref{avram}).

\begin{proposition}\label{K_22apr_1}
The value function $U$ of the McKean optimal stopping problem is given by

\[ U(x) = \left\{ \begin{array}{ll}
K-e^x & \mbox{if $x \leq k^*$} \\
c_1 e^{\beta_1(0,q)(x-k^*)} + c_2 e^{\beta_2(0,q)(x-k^*)} & \mbox{if $x > k^*$,}
\end{array} \right. \]
where

\begin{multline*}
c_1=\frac{\beta_2(0,q) K + (1-\beta_2(0,q)) e^{k^*}}{\beta_2(0,q)-\beta_1(0,q)}, \qquad c_2=\frac{\beta_1(0,q) K + (1-\beta_1(0,q)) e^{k^*}}{\beta_1(0,q)-\beta_2(0,q)} \\
\mbox{and} \qquad e^{k^*} = \frac{K q}{\sigma_{X}^2/2+q+\lambda/(\theta+1)^2}.
\end{multline*}

\end{proposition}

\begin{proof} A direct derivation of these formulas can be found in \cite{Kou02} e.g. Alternatively, plugging the formulas from Proposition \ref{K_7apr_1} in the results from Theorem \ref{avram} we see that we can write

\begin{equation}\label{K_1mei_1}
U(x) = K q \sum_{i=1}^3 \frac{C_i(0,q)}{\beta_i(0,q)}  e^{\beta_i(0,q) (x-k^*)} -e^x \quad \mbox{and} \quad e^{k^*} = K \frac{\psi(1)}{\psi'(1)}.
\end{equation}
Applying the identity

\begin{equation}\label{K_1mei_2}
\frac{\sigma_{X}^2}{2} \prod_{i=1}^3 (z-\beta_i(c,q)) = (\theta+z+c) (\psi_c(z)-q) \quad \mbox{for $z \not= -\theta-c$}
\end{equation}
to this particular case (i.e. $c=0$, $q=\psi(1)$, $\beta_3(0,q)=1$), dividing both sides by $z-1$ and taking the limit for $z \to 1$ we find

\begin{equation}\label{K_11mei1}
\sigma_{X}^2 (1-\beta_1(0,q)) (1-\beta_2(0,q)) = 2 (\theta+1) \psi'(1).
\end{equation}
Plugging this in the equation for $e^{k^*}$ we find $e^{k^*}=2(\theta+1)Kq/(\sigma_{X}^2 (\beta_2(0,q)-1)(\beta_1(0,q)-1))$. Using this expression in (\ref{K_1mei_1}), together with $\beta_1(0,q) \beta_2(0,q)=2q\theta/\sigma_{X}^2$ (from (\ref{K_1mei_2}) with $z=0$), the stated formula for $U$ indeed follows.
\end{proof}

Now we are ready to turn to formulas for the optimal exercise levels $x^*, y^*$ and the value function $V$ of the McKean game. Recall that for $\delta \geq U(\log K)$ the game degenerates to the McKean optimal stopping problem.

\begin{proposition}\label{K_1mei_3} Consider the McKean game driven by (\ref{K_6apr_1}). Recall $\bar \delta=U(\log K)$. We assume throughout that $\delta<\bar \delta$.
\begin{itemize}
\item[(i)] The optimal level $x^*=x^*(\delta)$ is the unique solution to the equation in $x$:

\[ q \sum_{i=1}^3 \frac{C_i(0,q)}{\beta_i(0,q)} K^{\beta_i(0,q)} e^{-\beta_i(0,q) x} -1 = \frac{\delta}{K}. \]
On $(-\infty,x^*]$ we have $V(x)=K-e^x$ and on $(x^*,\log K]$ we have

\[ V(x) = K q \sum_{i=1}^3 \frac{C_i(0,q)}{\beta_i(0,q)}  e^{\beta_i(0,q) (x-x^*)} -e^x. \]
\item[(ii)] The threshold $\delta_0 \in (0,\bar \delta)$ is the unique solution to the equation in $z$:

\[ q \sum_{i=1}^3 \frac{C_i(0,q) K^{\beta_i(0,q)}}{\beta_i(0,q) (\theta + \beta_i(0,q))} e^{-\beta_i(0,q) x^*(z)} - \frac{\lambda + (\theta+1)q}{\lambda \theta K} z = \frac{1}{\theta+1}. \]

\item[(iii)] Suppose $\delta \in [\delta_0,\bar \delta)$. We have $y^*=\log K$ and on $[\log K,\infty)$

\[ V(x) = K \sum_{i=1}^2  C_i(0,q) \left( \frac{q e^{-\beta_i(0,q) x^*}}{\beta_i(0,q)} + K^{-\beta_i(0,q)} \left( \psi'(1)-Kqe^{-x^*} \right) \right) e^{\beta_i(0,q) x}. \]

\item[(iv)] Suppose $\delta \in (0,\delta_0)$. We have

\[ e^{\theta y^*} = \frac{\lambda \theta K^{\theta+1}}{(\theta+1) q \delta} \left( q \sum_{i=1}^3 \frac{C_i(0,q) K^{\beta_i(0,q)}}{\beta_i(0,q) (\theta + \beta_i(0,q))} e^{-\beta_i(0,q) x^*} - \frac{1}{\theta+1} - \frac{\delta}{\theta K} \right). \]
On $[\log K,y^*]$ we have $V(x)=\delta$ and on $(y^*,\infty)$ we have

\[ V(x) = \frac{\delta}{\beta_2(0,q)-\beta_1(0,q)} \left( \beta_2(0,q) e^{\beta_1(0,q) (x-y^*)} - \beta_1(0,q) e^{\beta_2(0,q) (x-y^*)} \right). \]
\end{itemize}
\end{proposition}

\begin{proof} Ad (i). Apply Proposition \ref{K_7apr_1} to the formulas from Theorem \ref{mainthrm} (ii).

Ad (ii). Apply Proposition \ref{K_7apr_1} to Theorem \ref{thm_Erik2} (i).

Ad (iii). Apply Proposition \ref{K_7apr_1} to the formula from Theorem \ref{mainthrm} (ii) to obtain

\[ V(x) = K \sum_{i=1}^3  C_i(0,q) \left( \frac{q e^{-\beta_i(0,q) x^*}}{\beta_i(0,q)} + K^{-\beta_i(0,q)} \left( \psi'(1)-Kqe^{-x^*} \right) \right) e^{\beta_i(0,q) x}-e^x \]
and use (\ref{K_11mei1}) to see that the terms involving the exponential of a positive factor times $x$ vanish. (Of course, one can also reason directly that they should cancel, since otherwise $V$ would not stay bounded for large $x$, which it should by definition).

Ad (iv). For $y^*$, apply Proposition \ref{K_7apr_1} to Theorem \ref{thm_Erik2} (ii) and simplify to arrive at the stated formula. Note that

\[ \sum_{i=1}^3 \frac{C_i(0,q)}{\beta_i(0,q) (\theta + \beta_i(0,q))} = \frac{2}{\sigma_{X}^2 \prod_{i=1}^3 \beta_i(0,q)} = \frac{1}{\theta q}, \]
the final equality by (\ref{K_1mei_2}).

For $V$, apply Proposition \ref{K_7apr_1} to Theorem \ref{thm_Erik2} (ii) and simplify, making use of the formula for $y^*$ and in particular Proposition \ref{K_7apr_1} (ii).
\end{proof}

We conclude with some plots in this jump-diffusion setting to illustrate the main result from this paper.

Figures 1 and 2 show the value function $V$ in the two different cases $\delta \in [\delta_0,\bar{\delta})$ and $\delta \in (0,\delta_0)$.

Figure 3 shows a plot of $\bar{\delta}$ and $\delta_0$ as a function of $\sigma_X$. (Note that this really means only $\sigma_X$ changes, hence the equation $q=\psi(1)$ does not (necessarily) hold as $\psi$ changes with $\sigma_X$). This figure can be explained as follows. If $\sigma_X \downarrow 0$, $U$ converges to the value of the McKean optimal stopping problem for $X$ with $\sigma_X=0$ and hence $\bar{\delta}=U(\log K)$ also has some limit in the interval $(0,K)$. The figure suggests that the difference between $\bar{\delta}$ and $\delta_0$ vanishes as $\sigma_X \downarrow 0$, which might be explained as follows. As pointed out in Remark \ref{rem_Kees1}, $\delta \in (\delta_0,\bar{\delta})$ means that when $X$ starts above $\log K$, the probability of hitting $\log K$ before it reaches levels (far) below $\log K$ is large enough to have $y^*=\log K$. Obviously, this probability vanishes together with $\sigma_X$, and hence the length of the interval $(\delta_0,\bar{\delta})$ vanishes as $\sigma_X \downarrow 0$. Furthermore, as $\sigma_X \to \infty$, for the maximizer in the McKean optimal stopping problem the negative effect of discounting vanishes in the sense that first entry time of any interval has a density approaching the Dirac measure in $0$. Hence, $k^* \to -\infty$ and $U(x) \uparrow K$ for any $x$. In particular also $\bar{\delta}=U(\log K) \uparrow K$. The vanishing of $\bar{\delta}$ as $\sigma_X \to \infty$ is explained as above by the fact that increasing $\sigma_X$ means that the probability of hitting $\log K$ before $X$ falls (far) below $\log K$ increases, and hence the effect of the negative jumps is of vanishing relevance. In the limit therefore the minimizer does not need to choose an $y^*>\log K$ for any $\delta>0$.

Finally Figure 4 shows how $x^*$ and $y^*$ vary with $\sigma_X$. Note that the behaviour of $y^*$ is indeed consistent with the structure Figure 3 suggests.

\newpage

\begin{figure}[h!tb]
\begin{center}
\includegraphics{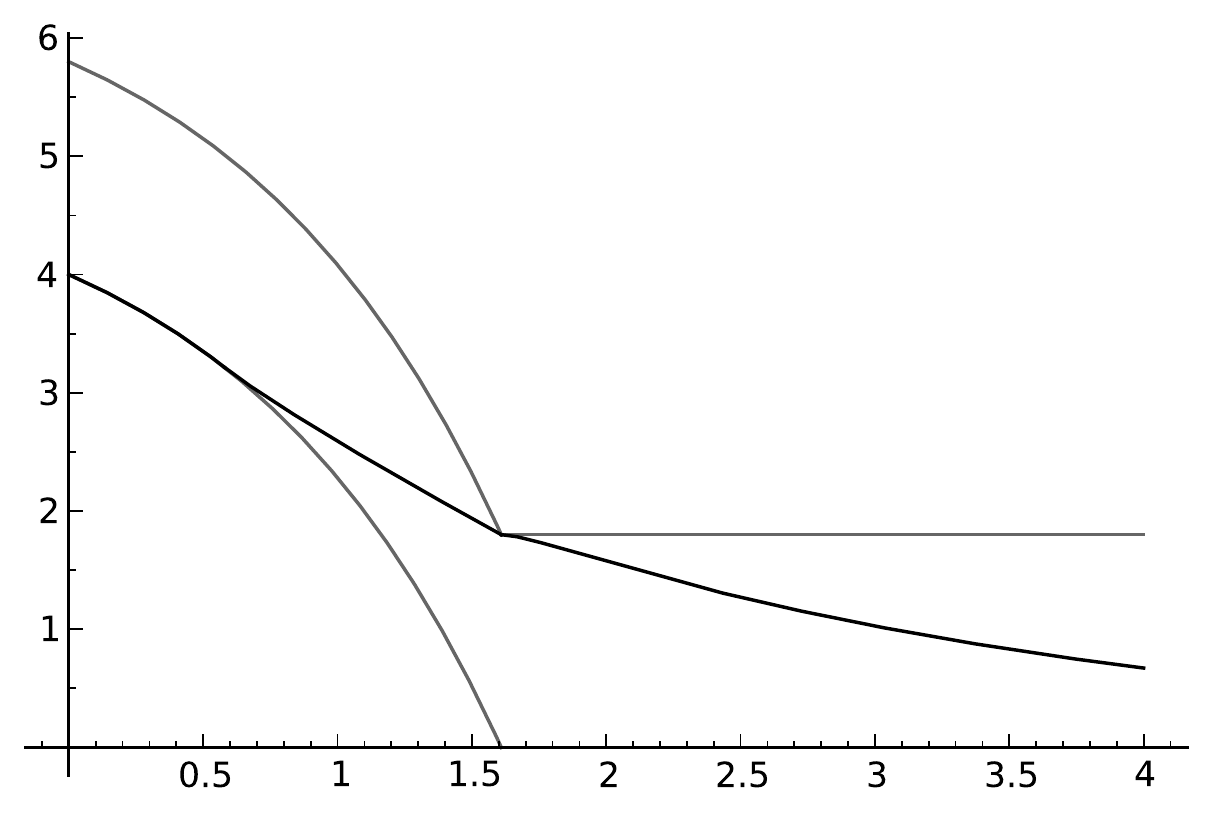}
\caption{A plot of the value function $V$ in the case $\delta \in [\delta_0,\bar \delta)$, so $y^*=\log K$. The grey curves are the upper and lower payoff functions, the black curve is $V$. Here $K=5$, $\delta=1.8$, $\bar{\delta}=2.08$, $\delta_0=1.76$ and $x^*=0.58$}
\end{center}
\end{figure}

\begin{figure}[h!tb]
\begin{center}
\includegraphics{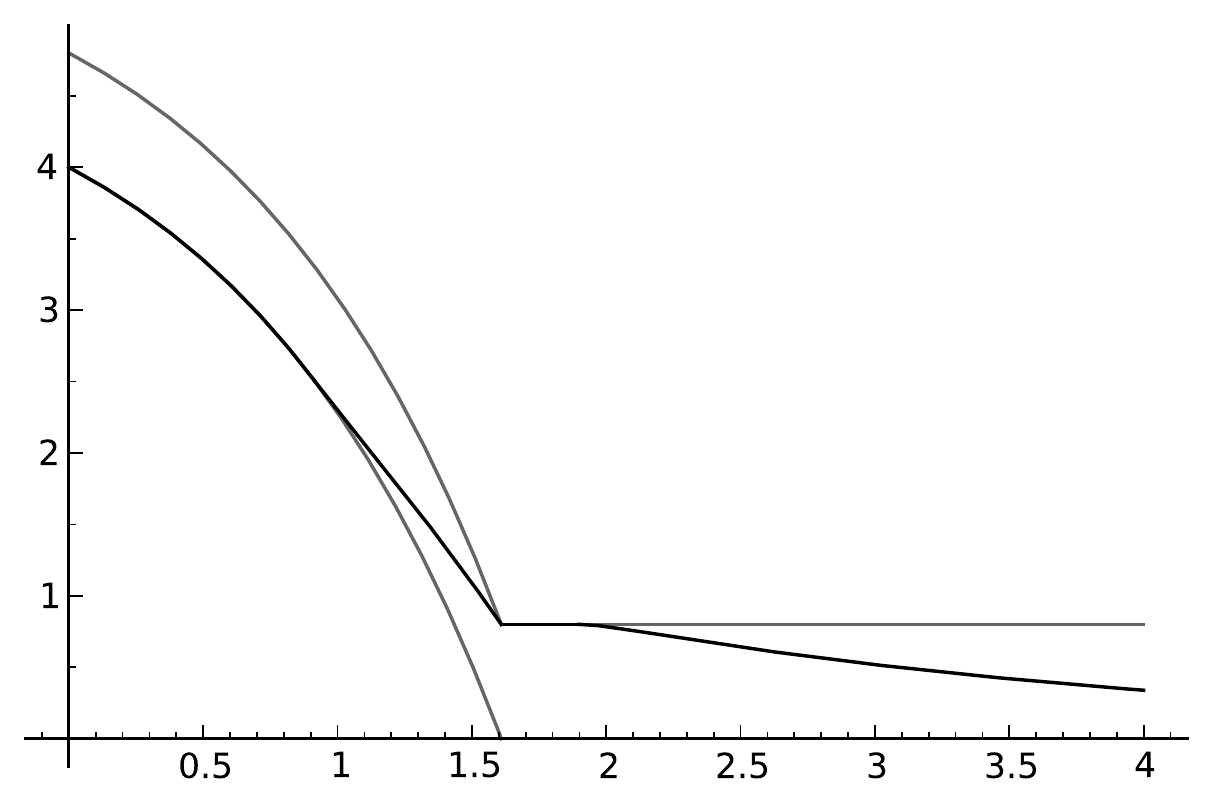}
\caption{A plot of the value function $V$ in the case $\delta \in (0,\delta_0)$, so $y^*>\log K$. The grey curves are the upper and lower payoff functions, the black curve is the value function $V$. Here $K$, $\delta_0$ and $\bar{\delta}$ are as in Figure 1; $\delta=0.8$, $x^*=0.91$ and $y^*=1.90$}
\end{center}
\end{figure}

\begin{figure}[h!tb]
\begin{center}
\includegraphics{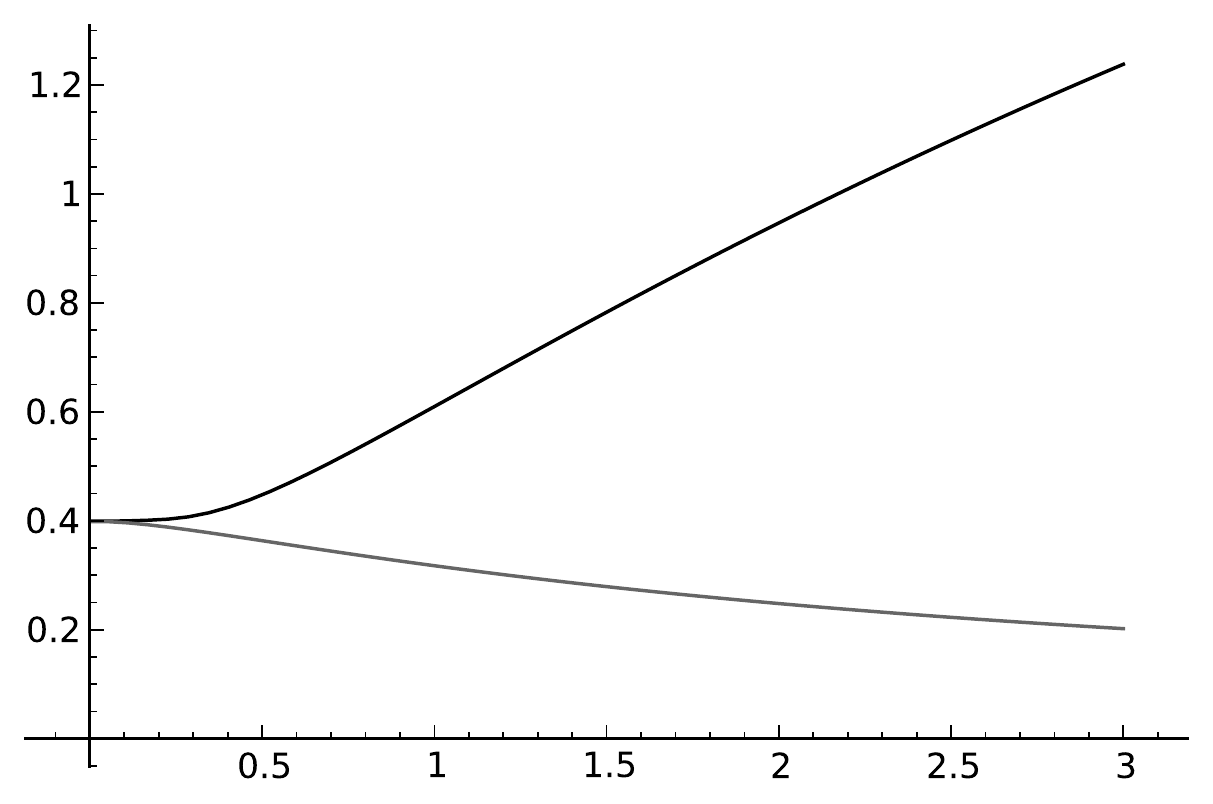}
\caption{A plot of $\bar{\delta}$ (the black curve) and $\delta_0$ (the grey curve) as a function of $\sigma_X$}
\end{center}
\end{figure}

\begin{figure}[h!tb]
\begin{center}
\includegraphics{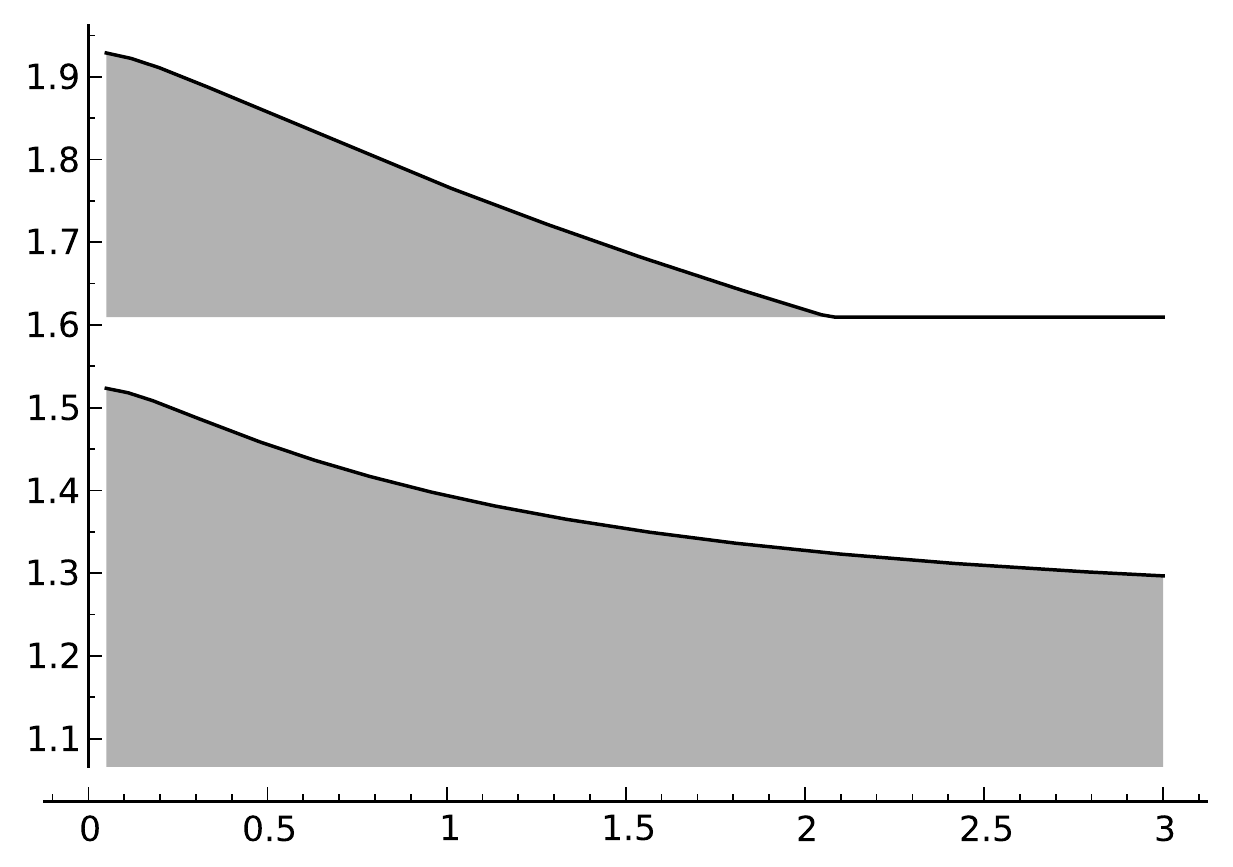}
\caption{A graphical illustration of the optimal stopping regions $(-\infty,x^*]$ (maximizer) and $[\log K,y^*]$ (minimizer) as a function of $\sigma_X$. The upper black curve is $y^*(\sigma_X)$, the lower one is $x^*(\sigma_X)$}
\end{center}
\end{figure}


\newpage

\begin{thebibliography}{99}
\bibitem{Avram04} Avram, F. and Kyprianou, A.E. and Pistorius, M. R. (2004) Exit problems for spectrally negative L\'evy processes and applications to (Canadized) Russian options. \textit{Ann. Appl. Probab.} \textbf{14}, 215--238.

\bibitem{McKean} Baurdoux, E.J. and Kyprianou, A.E. (2008) The McKean stochastic game driven by a spectrally negative L\'evy process. \textit{Elec. J. of Probab.} \textbf{8}, 173--197.

\bibitem{Baurdoux08} Baurdoux, E.J. and Kyprianou, A.E. (2008) The Shepp--Shiryaev stochastic game driven by a spectrally negative L\'evy process. \textit{To appear in Theory of Probability and Its Applications}.

\bibitem{Baurdoux09} Baurdoux, E.J. and Kyprianou, A.E. and Pardo, J.C. (2010) The Gapeev--K\"uhn stochastic game driven by a spectrally positive L\'evy process. \textit{Submitted}.

\bibitem{Bertoin96} Bertoin, J. (1996) L\'evy {P}rocesses. \textit{Cambridge University Press}.

\bibitem{Bichteler02} Bichteler, K. (2002) Stochastic Integration with jumps. \textit{Cambridge University Press} MR1906715.

\bibitem{Chan}  Chan, T. (2004) Some applications of L\'evy processes in
insurance and finance. \textit{Finance.}
{\bf 25}, 71--94.

\bibitem{Doney}
Doney, R.A. (2005) Some excursion calculations for spectrally one-sided L\'evy processes.
S\'em. Probab. \textbf{XXXVIII}, 5-15. 

\bibitem{Dynkin69} Dynkin, E. B. (1969) A game-theoretic version of an optimal stopping problem. \textit{Dokl. Akad. Nauk. SSSR} \textbf{185}, 16--19.

\bibitem{Ekstrom08} Ekstr\"om, E. and Peskir, G. (2008) Optimal stopping games for {M}arkov processes. \textit{SIAM J. Control Optim.} \textbf{2}, 684--702.

\bibitem{Gapeev05} Gapeev, P. V. and K\"uhn, C. (2005) Perpetual convertible bonds in jump-diffusion models. \textit{Statistics \& Decisions} \textbf{23}, 15--31

\bibitem{Kallsen04} Kallsen, J. and K\"uhn, C. (2004) Pricing derivatives of {A}merican and game type in incomplete markets. \textit{Finance and Stochastics} \textbf{8}, 261--284.

\bibitem{Kifer00} Kifer, Y. (2000) Game options. \textit{Finance and Stochastics} \textbf{4}, 443--463.

 \bibitem{Kou02}  Kou, S. and Wang, H. (2002) Option pricing under a jump-diffusion model. \textit{Management Science} \textbf{50}, 1178--1192.

 \bibitem{Kyprianou04} Kyprianou, A. E. (2004) Some calculations for {I}sraeli options. \textit{Finance and Stochastics} \textbf{8}, 73--86.

\bibitem{Kyprianou06} Kyprianou, A. E. (2006) Introductory Lectures on Fluctuations of L\'evy processes with Applications. \textit{Springer}.

\bibitem{McKean65} McKean, H. (1965) Appendix: A free boundary problem for the heat equation arising from a problem of mathematical economics. \textit{Ind. Manag. Rev.} \textbf{6}, 32--39.

\bibitem{Mordecki}  Mordecki, E. (2002) Optimal stopping and perpetual options for L\'{e}vy processes. \textit{Finance Stoch.} \textbf{6}, 473--493.

\end{thebibliography}
\end{document}